\documentclass{amsart}

\newtheorem{theorem}{Theorem}[section]
\newtheorem{proposition}[theorem]{Proposition}

\def\a{\alpha}
\def\b{\beta}

\def\F{\Phi}

\def\k{\kappa}

\def\z{\zeta}

\def\ov{\overline}

\def\wt{\widetilde}
\def\wh{\widehat}

\def\SZ{Szeg\H o}
\def\({\left(}
\def\){\right)}
\def\[{\left[}
\def\]{\right]}

   \def\dT{{\mathbb T}}

      \def\cC{{\mathcal C}}
      
\def\cG{{\mathcal G}}      
      \def\cL{{\mathcal L}}
\def\cM{{\mathcal M}}

\font\tenopen = cmbx10 \font\sevenopen = cmbx7 \font\fiveopen =
cmbx5
\newfam\openfam

\textfont\openfam = \tenopen \scriptfont\openfam = \sevenopen
\scriptscriptfont\openfam = \fiveopen

\begin{document}
\title{An inverse spectral theory for finite CMV matrices}

\author{Leonid Golinskii}
\address{Institute for Low Temperature Physics and Engineering\\47\\ Lenin
ave.\\ Kharkov, 61103 \\ Ukraine} \email{golinskii@ilt.kharkov.ua;
leonid.golinskii@gmail.com}

\author{Mikhail Kudryavtsev}
\address{Institute for Low Temperature Physics and Engineering\\47\\ Lenin
ave.\\ Kharkov, 61103 \\ Ukraine} \email{kudryavtsev@ilt.kharkov.ua;
kudryavstev@onet.com.ua}



\begin{abstract}
For finite dimensional CMV matrices the classical inverse spectral
problems are considered. We solve the inverse problem of
reconstructing a CMV matrix by its Weyl's function, the problem of
reconstructing the matrix by two spectra of CMV operators with
different ``boundary conditions'', and the problem of
reconstructing a CMV matrix by its spectrum and the spectrum of
the CMV matrix obtained from it by truncation. Bibliography : 24
references.
\end{abstract}

\maketitle

 \section{Introduction}\setcounter{section}{1}
\setcounter{equation}{0}

$ \ \  \ $     The theory of orthogonal polynomials on the unit
circle (OPUC) studies probability measures $\mu$ on the unit
circle $\mathbb{T} = \{|\zeta| =1 \} $, and polynomials which
emerge as an outcome of the Gram--Schmidt procedure applied to a
sequence of monomials $\{\zeta^j\}_{j\ge 0}$ in the Hilbert space
$L^2_\mu(\mathbb{T})$ with the inner product
$$ (f,g)_\mu = \int_{\dT} f(\zeta) \ov{g(\zeta)} d\mu \,, \qquad
\|f\|^2 = (f,f)_\mu \,.
$$
There are two natural ways of normalization: the orthonormal
polynomials
$$ \varphi_n(z) = \varphi_n(z,\mu) = {\k}_n z^n + \ldots, \quad
(\varphi_n,\varphi_m)_\mu=\delta_{n,m}, \qquad n,m \in
\mathbb{Z}_+=\{0,1,\ldots\},
$$
and the monic orthogonal polynomials
$$\Phi_n(z) = \Phi_n(z,\mu)={\k}_n^{-1}\varphi_n(z)+\ldots=z^n+\ldots,
\qquad (\Phi_n, \Phi_m) = 0, \ n\ne m \,.
$$Both systems are uniquely determined provided $\k_n >0$ is
required. For the background of the theory see \cite{Sz, Ger61,
simA, simB}.

An operator $U$ of multiplication by $\zeta$ in
$L^2_\mu(\mathbb{T})$, $Uf = \zeta f$, with $\mu$ being the
spectral measure associated to the constant function $f \equiv 1$,
is a backbone of the OPUC theory. But since $L^2_\mu(\mathbb{T})$
is $\mu$-dependent, such multiplication operators act in different
spaces, so one cannot connect them, which is crucial especially in
perturbation theory. So, a suitable matrix representation or, in
other words, a convenient orthonormal basis in $L^2_\mu(\dT)$ is
needed so that all the operators act ithe same space $\ell^2$.
There is an ``obvious'' set to try, namely, $\{\varphi_n\}$, but
it has two defects. First, by the fundamental \SZ\
--Kolmogorov--Krein theorem, $\{\varphi_n\}$ is a basis if and
only if $\mu$ is outside the \SZ\ class, that is, $\log \mu'
\not\in L^1(\mathbb{T})$, $\mu'$ is the Radon-Nikodym derivative
of $\mu$ with respect to the normalized Lebesque measure $dm$ on
$\mathbb{T}$. Second, even if this is the case, the corresponding
matrix $\cG$ is not of finite width measured from the diagonal.

One of the most interesting developments in the OPUC theory in
recent years is the discovery by Bunse-Gerstner and Elsner
\cite{BGE}, Watkins \cite{W93} and later (in a more transparent
form) by Cantero, Moral and Vel\'{a}zquez \cite{CMV2, CMV3, CMV4}
of a matrix realization for $U$ which is of finite band size.
Specifically, there exists an orthonormal basis $\{\chi_n\}$,
called the CMV basis, in $L^2_\mu(\mathbb{T})$ such that
$$ (\zeta \chi_m,\chi_n)_\mu = 0, \quad |m-n|>2,
$$
to be compared with Jacobi matrices which correspond to orthogonal
polynomials on the real line. $\{\chi_n\}$ is obtained by
orthonormalizing the sequence
$\{1,\zeta,\zeta^{-1},\zeta^2,\zeta^{-2},\ldots\}$, and the
mentioned above matrix realization
$$ \cC=\cC(\mu) = \|c_{n,m}\|_{n,m=0}^\infty, \qquad c_{n,m}=(\zeta
\chi_m,\chi_n)_\mu\,,
$$
called the {\it CMV matrix}, is unitary and five-diagonal.
Remarkably, the $\chi$'s can be expressed in terms of $\varphi$'s
and their reversed $\varphi^*$'s by
$$ \chi_{2n} = z^{-n} \varphi^*_{2n}(z), \quad
\chi_{2n+1}(z)=z^{-n}\varphi_{2n+1}(z), \qquad n\in \mathbb{Z}_+\,,
$$
where $\varphi_k^*(z) = z^k \overline{\varphi_k (1/\overline{z})}$.
\cite{CMV2} provides another crucial idea: $\cC$ can be factorized
into
\begin{equation}\label{1.1}
\cC=\cL \cM, \qquad \cL=\Theta_0 \oplus \Theta_2 \oplus \Theta_4
\oplus \ldots, \qquad \cM = 1 \oplus \Theta_1 \oplus \Theta _3
\oplus \ldots
\end{equation}
(the first block of $\cM$ is $1\times 1$), with
\begin{equation}\label{1.2}
\Theta_j=\begin{pmatrix} \overline{\alpha_j} & \rho_j \\
\rho_j & -\alpha_j \\
\end{pmatrix}, \quad j\in\mathbb{Z}_+\,, \qquad
\rho_j=\sqrt{1-|\a_j|^2}>0.
\end{equation}
The numbers $\alpha_j$ are known as the {\it Verblunsky
coefficients} of the measure $\mu_j$. By the Verblunsky theorem
(see \cite[Theorem 1.7.11]{simA}), $|\alpha_j|<1$ and, moreover,
each such sequence of complex number occurs.  Expanding out the
matrix product (\ref{1.1})--(\ref{1.2}), although rather
cumbersome, can be carried out and leads to a rigid structure
$$ \cC(\mu) = \begin{pmatrix}
* & * & + & & & & & & \\
+ & * & * & & & & & & \\
  & * & * & * & + & & & \\
  & + & * & * & * & & & \\
  &   &   & + & * & * & * & \\
\ldots & \ldots & \ldots & \ldots & \ldots & \ldots & \ldots &
\ldots\\
\end{pmatrix},
$$
where $+$ represents strictly positive entries, and $*$ generally
non-zero ones. The entries marked $+$ are precisely $(2,1)$ and
those of the form $(2j-1,2j+1)$ and $(2j+2,2j)$ with
$j\in\mathbb{N}=\{1,2,\ldots\}$, so the half of the entries
$c_{n,n+j}$ with $j=\pm2$ are zero, and $\cC$ is only "barely"
five-diagonal. The explicit formulae for $c_{m,n}$ in terms of
$\alpha$'s and $\rho$'s are also available (see \cite{Gol2}), in
particular,
\begin{equation}\label{1.3}
\cC(\mu) = \begin{pmatrix}
\bar\alpha_0 & \bar\alpha_1\rho_0 & \rho_1\rho_0 & 0 & 0 & \ldots  \\
\rho_0 & -\bar\alpha_1\alpha_0 & -\rho_1\alpha_0 & 0 & 0 & \ldots \\
0 & \bar\alpha_2\rho_1 & -\bar\alpha_2\alpha_1 & \bar\alpha_3\rho_2 & \rho_3\rho_2 & \ldots \\
0 & \rho_2\rho_1 & -\rho_2\alpha_1 & -\bar\alpha_3\alpha_2 & -\rho_3\alpha_2 & \ldots \\
0  &  0 & 0  & \bar\alpha_4\rho_3 & -\bar\alpha_4\alpha_3 & \ldots \\
\ldots & \ldots & \ldots & \ldots & \ldots & \ldots \\
\end{pmatrix}.
\end{equation}

While in the above construction $\alpha_n$'s come from $\mu$,
$\Theta_j$ (\ref{1.2}) define unitaries so long as $|\alpha_j|\le
1$. So, there is another point of view on the CMV matrices (see
\cite{sim2, KN}) as those of the form (\ref{1.3}) with arbitrary
$|\alpha_j|\le 1$. If $|\alpha_j|<1$ for all $j$, $\cC$ is called
a proper CMV matrix, and it is an object of OPUC theory.
Otherwise, if $|\alpha_j|=1$ for some $j$, $\cC$ is called an
improper CMV matrix (and it has no direct relation to OPUC
theory). Interestingly enough, when
$|\alpha_0|,\ldots,|\alpha_{n-2}|<1$, and $|\alpha_{n-1}|=1$, we
have
\begin{equation}\label{1.4}
\cC=\cC_n\oplus \cC_\infty\,,
\end{equation}
where $\cC_n= \cC(\alpha_0,\ldots,\alpha_{n-2};\alpha_{n-1})$ is
$n\times n$ unitary matrix called a {\it finite CMV matrix}. The
class of such matrices, parameterized by an arbitrary set
$(\alpha_0,\ldots,\alpha_{n-2};\beta)$ with $|\alpha_j|<1$,
$|\beta|=1$, is the main item of business of the present paper.
For instance, if in (\ref{1.3}) $|\alpha_3|=1$, so $\rho_3=0$, we
have a finite CMV matrix of order $4$ as a principal $4\times 4$
block.

There is a multiplication formula for finite CMV matrices, similar
to (\ref{1.1}), which now depends on the parity of $n$:
$$ \cC(\alpha_0,\ldots,\alpha_{2k-1};\beta)=\cL_{2k+1}\cM_{2k+1}$$
with
$$\cL_{2k+1} = \Theta(\alpha_0) \oplus \ldots \oplus
\Theta(\a_{2k-2})\oplus\bar\b, \quad \cM_{2k+1} =
1\oplus\Theta(\alpha_1) \oplus \ldots \oplus \Theta(\a_{2k-1}), $$
and
$$ \cC(\alpha_0,\ldots,\alpha_{2k-2};\beta)=\cL_{2k}\cM_{2k}$$
with
$$\cL_{2k} = \Theta(\alpha_0) \oplus \ldots \oplus
\Theta(\a_{2k-2}), \quad \cM_{2k+1} = 1\oplus\Theta(\alpha_1)
\oplus \ldots \oplus \Theta(\a_{2k-1})\oplus\bar\b. $$

Our argument is based upon two results originated in OPUC theory.
The first one is the famous {\it \SZ\ recurrence relations}
\begin{equation}
\label{1.5} \Phi_k(z)=z\Phi_{k-1}(z) - \bar\alpha_{k-1}
\Phi_{k-1}^*(z); \qquad k=1,2,\ldots,n-1, \quad \Phi_0\equiv 1,
\end{equation}
where $\Phi_k = \det (z-\cC^{(n)})$, $\cC^{(n)}$ a principal
$k\times k$ submatrix of a finite CMV matrix
$\cC=\cC(\alpha_0,\ldots,\alpha_{n-2};\beta)$. We single out the
final relation:
\begin{equation}
\label{1.6} \wt\Phi_n(z)=z\Phi_{n-1}(z) - \bar{\beta}
\Phi_{n-1}^*(z); \quad \wt\Phi_n = \det(z-\cC),
\end{equation}
to emphasize that $\wt\Phi_n$ is no longer an ``orthogonal'', but
a ``paraorthogonal'' polynomial:
\begin{equation}
\label{1.7} \wt\Phi_n(z) = \prod_{j=1}^n (z-\zeta_j), \quad
\zeta_j \in \mathbb{T}.
\end{equation}
Note that one can prove (\ref{1.5})--(\ref{1.6}) by direct
expanding out the determinants (no orthogonality is needed here).
In what follows we call $\Phi_k$, $\wt\Phi_n$ the {\it \SZ\
polynomials} associated with a finite CMV matrix $\cC$.

The second result, known as Geronimus' theorem (see \cite[Theorem
1.7.5]{simA}), reads that given a monic polynomial $P_k$ of degree
$k$ with all its zeroes inside the unit disk
$\mathbb{D}=\{|z|<1\}$ (as it belongs for OPUC), there is a
measure $\mu$ such that $P_k=\Phi_k(\mu)$. What is more to the
point, $\Phi_k(\mu_1)=\Phi_k(\mu_2)$ implies
$\alpha_j(\mu_1)=\alpha_j(\mu_2)$ for $j=0,1,\ldots,k-1$, and
$\alpha_j$ can be reconstructed from the {\it inverse \SZ\
recurrence} \cite[Theorem 1.5.4]{simA}
\begin{equation}\label{invszeg}
z\Phi_j(z)=\rho_j^{-2}\left(\Phi_{j+1}(z)+\bar\alpha_j
\Phi_{j+1}^*(z)\right), \quad \a_j=-\ov{\F_{j+1}(0)}, \quad
j=0,1,\ldots,n-2,
\end{equation}
In particular, for $\cC=\cC(\alpha_0,\ldots,\alpha_{n-2};\beta)$
the polynomial $\Phi_{n-1}$ completely determines
$\alpha_0,\ldots,\alpha_{n-2}$, and $\wt\Phi_n$ produces $\beta$
from (\ref{1.6}) and (\ref{1.6})--(\ref{1.7}) with $z=0$,
$(\F_{n-1}^*(0)=1)$:
\begin{equation}
\label{1.8} \beta=-\overline{\wt\Phi_n(0)} = (-1)^{n+1}
\prod_{j=1}^n \bar{\zeta_j}.
\end{equation}
Hence, $\cC$ is uniquely determined by two polynomials
$\Phi_{n-1}$ and $\wt\Phi_n$.

It is known \cite{sim2, KN} that each finite CMV matrix has a
simple spectrum, that is, all its eigenvalues $\{\zeta_j\}_1^n$
are distinct: $\zeta_j=\zeta_k$, $j\ne k$, and, moreover, each
$n\times n$ unitary matrix with a simple spectrum is unitary
equivalent to some $n\times n$ CMV matrix. By the Spectral Theorem
for unitaries
$$ \cC = \sum_{j=1}^n \zeta_j\Pi_j = \sum_{j=1}^n \zeta_j(\cdot,
h_j)h_j\,,
$$
where $\{h_j\}_1^n$ is the orthonormal basis in $\mathbb{C}^n$ of
the eigenvectors of $\cC$. Denote by $\{e_j\}_1^n$ the standard
basis in $\mathbb{C}^n$. The {\it Weyl function} ($n$-th Weyl
function) is then defined by
\begin{equation} \label{1.9}
w(\cC,z)=w(z) = \Bigl( (z-\cC)^{-1} e_n,e_n\Bigr) = \sum_{j=1}^n
{\mu_j \over z-\zeta_j}\,, \quad \mu_j=|(h_j,e_n)|^2.
\end{equation}
By the Cramer rule
\begin{equation} \label{1.10}
w(z) = {\Phi_{n-1}(z) \over \wt\Phi_n(z)} = \sum_{j=1}^n {1 \over
z-\zeta_j} \,  {\Phi_{n-1}(\zeta_j) \over \wt\Phi_n'(\zeta_j)}
\end{equation}
and since $\Phi_{n-1}(\zeta_j) \ne 0$ for all $j$,
\begin{equation} \label{1.11}
\mu_j=|(h_j,e_n)|^2 = {\Phi_{n-1}(\zeta_j) \over
\wt\Phi_n'(\zeta_j)}
>0, \qquad  \sum_{j=1}^n \mu_j = \|e_j\|^2 =1.
\end{equation}
We will also need the {\it $n$-th spectral measure} of $\cC$, that
is, associated with $e_n$ and defined by
\begin{equation} \label{1.12}
d\mu(\cC) = d\mu = \sum_{j=1}^n \mu_j \delta(\zeta_j), \qquad w(z)
= \int_\mathbb{T} {d\mu(\zeta) \over z-\zeta}.
\end{equation}

The inverse problems considered in the presented work are well
known for the Jacobi matrices. We proceed as follows. In Section 2
we consider the inverse problem of reconstructing the CMV matrix
from its Weyl function and from the spectra of two CMV matrices
with different ``boundary conditions''. In Section 3 we study the
spectra of CMV matrices and their truncations and solve the direct
problem, i.e., we find the necessary properties for two spectra to
be the data of the corresponding inverse problem. Finally, in
Section 4 we solve the inverse problem of reconstructing the CMV
matrix by the two spectra, examined in the preceding section. In
our forthcoming paper \cite{GK2} we consider the mixed inverse
problems for finite CMV matrices.

Our study of the inverse spectral problems for finite CMV matrices
was strongly inspired and influenced by the similar theory for
Jacobi matrices, developed in various papers (see, e.g.
\cite{Boley-Golub, Boor-Golub, Gray-Wilson, Hald, Hochstadt1,
Hochstadt2} and the survey \cite{GS} for the modern approach) as
well as the recent book \cite{Mar} for physically motivated
exposure of the inverse theory for Jacobi matrices.

\section{Inverse problems by spectral measure and two spectra}
\setcounter{section}{2} \setcounter{equation}{0}

As a warmup, we begin with the $n$-th spectral measure
(\ref{1.12}).

Given $\cC=\cC(\alpha_0,\ldots,\alpha_{n-2};\beta)$ we find its
\SZ\ polynomials $\{\Phi_k\}_1^{n-1}$ and $\wt\Phi_n$ by the \SZ\
recurrences (\ref{1.5})--(\ref{1.6}). Next, the zeroes
$\{\zeta_j\}_1^n$ of $\wt\Phi_n$ are the eigenvalues of $\cC$, and
the residues (\ref{1.11}) of the Weyl function $w$ are the masses,
so the direct spectral problem of finding the $n$-th spectral
measure from the CMV matrix is solved. Conversely, starting out
from the spectral measure $\mu$ (or, equivalently, from the Weyl
function $w$), we write
$$ \wt\Phi_n(z) = \prod_{j=1}^n(z-\zeta_j), \quad
\Phi_{n-1}(z)= \wt\Phi_n(z) w(z) = \wt\Phi_n(z) \sum_{j=1}^n
{\mu_j \over z - \zeta_j}\,,
$$ and compute $\alpha_0,\ldots,\alpha_{n-2}; \beta$ from
(\ref{invszeg}) and (\ref{1.8}). So $\cC$ can be easily retrieved
from its $n$-th spectral measure.

\begin{theorem}\label{theor1}
Given a set $\{\xi_j\}_1^n$ of $n$ distinct points on
$\mathbb{T}$, and a set $\{\nu_j\}_1^n$ of positive numbers with
$\sum_{j=1}^n \nu_j =1$, there is a unique CMV matrix $\cC$ so
that its $n$-th spectral measure is
$$d\mu(\cC) = \sum_1^n \nu_j \delta(\xi_j).
$$
\end{theorem}

{\it Proof.} \ Define
$$ P_n(z) = \prod_{j=1}^n (z-\xi_j), \quad \omega(z)=\sum_{j=1}^n {\nu_j
\over z-\xi_j }\,,
$$
and put
\begin{equation} \label{2.0}
P_{n-1}(z)=\omega(z) P_n(z) = ( \sum_{j=1}^n \nu_j ) z^{n-1}+\ldots
= z^{n-1}+\ldots.
\end{equation}

We show first that all zeroes of $\omega$ (equivalently, of
$P_{n-1}$) are in $\mathbb{D}$. In fact, let
$$ \omega(\xi_0)=\sum_1^n {\nu_j \over \xi_0-\xi_j } =0.
$$
Clearly, $\xi_0\ne\xi_j$. Then
$$\overline{\omega(\xi_0)} = \sum_1^n {\nu_j \over \bar\xi_0 -
\bar\xi_j } = \sum_{j=1}^n {\nu_j (\xi_0 - \xi_j) \over
|\xi_0-\xi_j|^2 } = 0,
$$
and so,
$$\xi_0 \sum_{j=1}^n {\nu_j \over |\xi_0-\xi_j|^2} =
\sum_{j=1}^n \xi_j {\nu_j \over |\xi_0-\xi_j|^2} ; \qquad \xi_0
=\sum_{j=1}^n\sigma_j\xi_j
$$
with
$$\sigma_j= {\nu_j \over
{|\xi_0-\xi_j|^2 {\displaystyle \sum_{j=1}^n}\nu_j
|\xi_0-\xi_j|^{-2}}}
> 0, \qquad \sum_{j=1}^n \sigma_j =1.
$$
The latter means that $\xi_0$ belongs to the convex hall of
$\{\xi_j\}_1^n$, and therefore $\xi_0\in\mathbb{D}$, as claimed.

By Geronimus' theorem there is a unique set of parameters
$\alpha_0,\ldots,\alpha_{n-2} \in \mathbb{D}$, and a system of the
\SZ\ polynomials $\{\Phi_0,\ldots,\Phi_{n-1}\}$ such that
(\ref{1.5}) holds, and $P_{n-1}=\Phi_{n-1}$. Put
\begin{equation}\label{2.00}
\wt\Phi_n(z) := z P_{n-1}(z) - \bar\beta P_{n-1}^*(z), \qquad
\beta = (-1)^{n+1} \prod_{j=1}^n \bar\xi_j. \end{equation}
 The CMV matrix $\cC_l =
 \cC(\alpha_0,\ldots,\alpha_{n-2};\beta)$ is now determined. Let
us show that $\omega$ is its Weyl's function $w(\cC,z)$.

We have by (\ref{2.0})
$$ \begin{aligned} z P_{n-1}(z) &= z \omega(z)
P_n(z) = \sum_{j=1}^n {z\nu_j \over z-\xi_j} P_n(z), \\
P_{n-1}^* &= z^{n-1} \overline{P_{n-1}(1/\overline{z})} = z^{-1}
\overline{\omega(1/\overline{z})}P_n^*(z) = \omega_1(z) P_n^*(z)
\end{aligned} $$
with
$$ \omega_1(z) = {1\over z} \sum_{j=1}^n{\nu_j \over z^{-1} - \xi_j} =
\sum_{j=1}^n {\nu_j\xi_j \over \xi_j-z} .
$$
Since ${\displaystyle P_n^*(z) = z^n \prod_1^n (z^{-1}
-\overline{\xi_j}) = (-1)^n \prod_1^n \bar\xi_j\, P_n(z)}$, then
by (\ref{2.00})
$$\wt\Phi_n(z) = \Bigl\{ \sum_1^n {z\nu_j\over z-\xi_j} -
\bar\beta \sum_1^n {\nu_j\xi_j \over \xi_j-z}\, (-1)^n \prod_1^n
\bar\xi_j \Bigr\} P_n(z) = P_n\,.
$$
Finally (see (\ref{1.10}))
$$ \omega(z) = {P_{n-1}(z)\over P_n(z)} = {\Phi_{n-1}(z)\over
\wt\Phi_n(z)} = w(\cC,z),$$ which completes the proof.  \hfill
$\square$

\bigskip

Let us go over to the inverse spectral problem related to two CMV
matrices with distinct ``boundary conditions''
\begin{equation} \label{2.1}
\cC_l = \cC(\alpha_0,\ldots,\alpha_{n-2};\beta_l), \quad l=1,2,
\quad \beta_1\ne\beta_2\,.
\end{equation}
It was proved in \cite{CMV1, Gol1} that the spectra $\Sigma(\cC_l)
= \{\zeta_{j,l}\}_{j=1}^n$ of such matrices interlace (for the
precise definition of intermittency of two point sets on
$\mathbb{T}$, which is quite consistent with the intuition, see
\cite{sim3}). Here is a simple reasoning. By (\ref{1.6})
$$ b(\zeta_{j,l})=\bar\beta_l,
\qquad j=1,2,,\ldots,n, \quad l=1,2; \qquad
b(z)={z\Phi_{n-1}(z)\over\Phi_{n-1}^*(z)}
$$
is a finite Blaschke product. As it is known, for
$\tau_1\ne\tau_2$ on $\mathbb{T}$, the $\tau_1$--points (the
points where $b(z)=\tau_1$) of any finite Blaschke product
interlace with its $\tau_2$--points, as needed.

Conversely, let the spectra $\Sigma(\cC_l)$ of (\ref{2.1}) be
known. Then by (\ref{1.8}) \newline
$\beta_l=(-1)^{n+1}\prod_j\bar\zeta_{j,l}$ are determined. Next,
by (1.6)-(1.7)
$$\wt\Phi_{n,1}(z)-\wt\Phi_{n,2}(z) = \prod_{j=1}^n(z-\zeta_{j,1})
- \prod_{j=1}^n (z-\zeta_{j,2}) =
(\bar\beta_2-\bar\beta_1)\Phi_{n-1}^*,
$$
so we find $\Phi_{n-1}$, and by Geronimus' theorem restore the
rest of Verblunsky parameters $\alpha_0,\ldots,\alpha_{n-2}$.

\begin{theorem}\label{theor2}
Let $\{\xi_{j,l}\}_{j=1}^n$, $l=1,2$, be two interlacing $n$-point
sets on $\mathbb{T}$. There is a unique pair of CMV matrices
$(\ref{2.1})$ such that
$$\Sigma(\cC_l)=\{\xi_{j,l}\}_{j=1}^n\,, \quad l=1,2.
$$
\end{theorem}

{\it Proof.} Put
 $$ P_{n,l}(z)=\prod_{j=1}^k(z-\xi_{j.k}), \quad l=1,2. $$
By \cite[Theorem 11.5.6]{simB} the difference $P=P_{n,1}-P_{n,2}$
has no zeroes in the closed unit disk. Since $\deg P \leq n-1$ and
$P(0)\ne 0$, then
$$ P_{n-1}(z):= {P^*(z)\over \overline{P(0)}}=z^{n-1}+\ldots
$$
is a monic polynomial of degree $n-1$ with all its zeroes in
$\mathbb{D}$. By Geronimus' theorem the set of parameters
$\alpha_0,\ldots,\ldots\alpha_{n-2}$ and the system of \SZ\
polynomials $\{\F_0,\ldots,\F_{n-1}\}$ can be reconstructed, and
$\F_{n-1}=P_{n-1}$. If we put $\beta_l = (-1)^{n+1} \prod_j
\bar\xi_{j,l}$, we end up with CMV matrices $\cC_l$ (\ref{2.1}).

It remains only to show that
$\Sigma(\cC_l)=\{\xi_{j,l}\}_{j=1}^n$, $l=1,2$. But
$\Sigma(\cC_k)$ are zero sets of polynomials
$$ \wt\Phi_{n,l}(z) = z P_{n-1}(z) - \bar\beta_l
P_{n-1}^*(z), \quad l=1,2.
$$
We have now
$$ P_{n-1}(z)={1\over \overline{P(0)}} z^{n-1}
\bigl\{ \overline{P_{n,1}(1/\overline{z})} -
\overline{P_{n,2}(1/\overline{z})} \bigr\}, \qquad \bar\beta_l = -
P_{n,l}(0) \in \mathbb{T},
$$
and $ z^{n-1} \overline{ P_{n,l} (1/\overline{z}) } = z^{-1}
\overline{P_{n,l}(0)} P_{n,l}(z)$, so
$$ \begin{aligned} z P_{n-1}(z) &=
{1\over \overline{P(0))}} \bigl\{ \overline{P_{n,1}(0)} P_{n,1}(z)
- \overline{P_{n,2}(0)} P_{n,2}(z) \bigr\} \\
&={P_{n,1}(0)P_{n,2}(0)\over P_{n,2}(0)-P_{n,1}(0)} \Bigl\{
{P_{n,1}(z)\over P_{n,1}(0)} - {P_{n,2}(z)\over P_{n,2}(0)}
\Bigr\} \\ &= {P_{n,2}(0)P_{n,1}(z) - P_{n,1}(0)P_{n,2}(z) \over
P_{n,2}(0)-P_{n,1}(0)}. \end{aligned}
$$
Next,
$$\bar\beta_l P_{n-1}^*(z) = - P_{n,l}(0){P(z)\over P(0)}
= P_{n,l}(0) {P_{n,1}(z) - P_{n,2}(z) \over P_{n,2}(0) -
P_{n,1}(0)},
$$
so, finally $\wt\Phi_{n,l}=P_{n,l}(z)$, as needed. \hfill
$\square$

\medskip

{\bf Example 1}. \ Let $\cC(\alpha_0,\ldots,\alpha_{n-2};\beta)$
be a finite CMV matrix with the \SZ\ polynomials
$\{\Phi_0,\ldots,\Phi_{n-1};\wt\Phi_n\}$, and let
$\tau\in\mathbb{T}$. Consider a new system \newline
$\{\Phi_0(\tau),\ldots,\Phi_{n-1}(\tau);\wt\Phi_n(\tau)\}$:
$$ \begin{aligned} \Phi_k(z,t)&=\tau^{-k}\Phi_k(\tau z),\quad
k=0,\ldots,n-1, \\
\wt\Phi_n(z,\tau)&=\tau^{-n}\wt\Phi_n(\tau z) = \prod_{j=1}^n
(z-\tau^{-1}\zeta_j). \end{aligned}
$$
It is clear that $\Phi_k(\tau)$ satisfy (\ref{1.5})
$$\Phi_k(z,\tau) = z \Phi_{k-1} (z,\tau) -
\overline{\alpha_{k-1}(\tau)}\Phi_{k-1}^*(z,\tau), \quad
\wt\Phi_n(z,\tau) = z \Phi_{n-1} (z,\tau)-
\overline{\beta(\tau)}\Phi_{n-1}^*(z,\tau)
$$
with
$$\alpha_k(\tau)=\alpha_k \tau^{k+1}, \quad k=0,\ldots,n-2; \quad
\beta(\tau)=\beta\tau^n. $$ Since
$\zeta_j(\tau)=\tau^{-1}\zeta_j$, $j=1,\ldots,n$, we have the
rotation of the spectrum.

When $\alpha_0=\ldots=\alpha_{n-2}=0, \ \beta=1$, then
$\wt\Phi_n(z)=z^n-1$, and for $\cC_0=\cC(0,\ldots,0;1)$ the
spectrum $\Sigma(\cC_0)=\bigl\{e^{2\pi ij\over n}\bigr\}_{j=1}^n$
consists of the $n$-th roots of unity. For $0<t<{2\pi\over n}$,
$\tau=e^{it}$, and $\cC_\tau=\cC(0,\ldots,0;\tau^n)$ we have the
rotation of that set, and for $\tau_1\ne\tau_2$ the spectra
$\Sigma(\cC_{\tau_1}), \Sigma(\cC_{\tau_2})$ interlace.

\section{Truncated CMV matrices, direct problem}
\setcounter{section}{3} \setcounter{equation}{0}

We examine here the CMV analog of the known inverse problem of the
reconstruction of an $n\times n$ Jacobi matrix from its spectrum and
the spectrum of its principal submatrix of order $n-1$. Consider two
CMV matrices
\begin{equation} \label{3.1}
\cC_1 = \cC(\alpha_0,\ldots,\alpha_{n-3},\alpha_{n-2};\beta_1),
\qquad \cC_2 = \cC(\alpha_0,\ldots,\alpha_{n-3};\beta_2),
\end{equation}
of order $n$ and $n-1$, respectively, and call $\cC_2$ a
"truncation" of $\cC_1$ (so, given $\cC_1$, there is a family of
truncations depending on a parameter $\beta_2\in\mathbb{T}$).

Let $\Sigma(\cC_1)=\{ \zeta_{j,1} \}_{j=1}^n$ and
$\Sigma(\cC_2)=\{\zeta_{j,2}\}_{j=1}^{n-1}$ be their spectra. The
direct spectral problem for the pair of matrices $\cC_l$
(\ref{3.1}) is solved in \cite{sim3}. We suggest an alternative
proof of this result, and develop some technique, which is applied
for the solution of the corresponding inverse spectral problem.

Given two $m$-point sets $Z_l=\{z_{j,l}\}_{j=1}^m$, $l=1,2$, on
$\mathbb{T}$, we will label them in the order
\begin{equation} \label{3.2}
z_{j,l} =e^{ix_{j,l}}: \quad 0\leq x_{1,l} < x_{2,l} <
\ldots<x_{m,l}<2\pi
\end{equation}
If $Z_1$ and $Z_2$ have no common points, we can assume that
$x_{1,1}<x_{1,2}$. So $Z_1$ and $Z_2$ interlace if
\begin{equation} \label{3.3}
0 \leq x_{1,1} < x_{1,2} < x_{2,1} < x_{2,2} < \ldots < x_{m,1}<
x_{m,2}\,.
\end{equation}
The following simple characterization of interlacing is crucial for
the rest of the paper. To simplify the notations, let us put
$$ \prod^{n(p)} a_j := a_1 a_2 \ldots a_{p-1} a_{p+1} \ldots a_n
=a_1a_2\ldots\wh a_p\ldots a_n\,.
$$
\begin{proposition} \label{propos3}
Let $Z_l=\{z_{j,l}\}_{j=1}^m$, $l=1,2$, be two point sets on
$\mathbb{T}$ with no common points, labeled by $(\ref{3.2})$. Put

\begin{equation} \label{3.4}
\omega_k := \frac{\prod\limits_{j=1}^m
\displaystyle\sin\frac{x_{k,1}-x_{j,2}}{2}} {\prod\limits_{}^{m(k)}
\displaystyle\sin\frac{x_{k,1}-x_{j,1}}{2}} = {u_k \over v_k}\,,
\qquad \omega_k\ne0, \quad k=1,2,\ldots,m.
\end{equation}
$Z_1$ and $Z_2$ interlace if and only if all $\omega_k$ have the
same sign.
\end{proposition}

{\it Proof.} \ With no loss of generality we can assume
$x_{1,1}<x_{1,2}$. It is easily seen from (\ref{3.3}) that
interlacing implies $\omega_k<0$ for all $k$. Conversely, let
$\omega_k<0$, $k=1,2,\ldots,m$. Define
$$ p_k:= \#\{ x_{j,2}: \ x_{k,1}<x_{j,2}<x_{k+1,1}\}; \qquad
k=1,2,\ldots,m, \quad x_{m+1,1}=x_{1,1}+2\pi,
$$
to be the number of points of $Z_2$ between $x_{k,1}$ and
$x_{k+1,1}$. Clearly, $p_k\in \mathbb{Z}_+$ and $\sum_k p_k =m$.
From the definition (\ref{3.4}) of $\omega_k$ one has
$$ {\rm sgn} \, v_{m-k} = (-1)^k, \quad {\rm sgn} \, u_{m-k} =
(-1)^{p_{m-k}+\ldots+p_m}.
$$
Since $\omega_k<0$ then $p_m$ is an odd number, $p_m+p_{m+1}$ is
an even number, $p_m+p_{m-1}+p_{m-2}$ is an odd number again, and
so on. Hence $p_j>0$, and $\sum_1^m p_k=m$ gives $p_k=1$ for all
$k$, that is exactly interlacing. \hfill $\square$

\medskip

The direct spectral result below is due to Simon \cite{sim3}.

\begin{theorem}\label{theor4}
Let $\cC_1$ and $\cC_2$ be two CMV matrices $(\ref{3.1})$, with
the spectra $\Sigma(\cC_1)=\{\zeta_{j,1}\}_{j=1}^n$ and
$\Sigma(\cC_2)=\{\zeta_{j,2}\}_{j=1}^{n-1}$. Define a unimodular
complex number

\begin{equation} \label{3.5}
B:= \bar\beta_1\bar\beta_2 {\beta_2-\alpha_{n-2} \over \bar\beta_2
- \bar\alpha_{n-2} } = \bar\beta_1 \beta_2 {1 -
\bar\beta_2\alpha_{n-2} \over 1 - \beta_2 \bar\alpha_{n-2} } \in
\mathbb{T}.
\end{equation}
Then the following dichotomy holds:\\
{\rm (i)} $\Sigma(\cC_1)\cap\Sigma(\cC_2)=\emptyset$, and the
$n$-point sets
$Z_1 = \Sigma(\cC_1)$, $Z_2=\Sigma(\cC_2)\cup\{B\}$ interlace;\\
{\rm (ii)} $\Sigma(\cC_1)\cap\Sigma(\cC_2)=\{ B\}$ and the
$(n-1)$-point sets $Z_1=\Sigma(\cC_1)\backslash\{B\}$,
$Z_2=\Sigma(\cC_2)$ interlace.
\end{theorem}

{\it Proof.} \ Let
$\{\Phi_0,\Phi_1,\ldots,\Phi_{n-2},\Phi_{n-1};\wt\Phi_n\}$ and
$\{\Phi_0,\Phi_1,\ldots,\Phi_{n-2};\wt\Phi_{n-1}\}$ be two systems
of the \SZ\ polynomials associated with $\cC_1$ and $\cC_2$,
respectively. Write the \SZ\ recurrences (\ref{1.5})--(\ref{1.6})
for $\Phi_{n-1}$ and $\wt\Phi_{n-1}$:
$$
\left\{\begin{array}{l}
\Phi_{n-1}(z) = z \Phi_{n-2}(z) - \bar\alpha_{n-2} \Phi_{n-2}^* (z),\\
\wt\Phi_{n-1}(z) = z \Phi_{n-2}(z) - \bar\beta_2 \Phi_{n-2}^* (z),\\
\end{array}
\right.
$$
and so
\begin{equation} \label{3.6}
\Phi_{n-1}(z) = \wt\Phi_{n-1}(z) + ( \bar\beta_2 -
\bar\alpha_{n-2}) \Phi_{n-2}^*(z).
\end{equation}
Next, write (\ref{1.6}) for $\wt\Phi_n$:
$$ \begin{aligned} \wt\Phi_n &= z\Phi_{n-1}(z) -
\bar\beta_1\Phi_{n-1}^*(z)\\
&= z(z\Phi_{n-2}(z) - \bar\alpha_{n-2} \Phi_{n-2}^*(z))
-\bar\beta_1 (\Phi_{n-2}^*(z) - z \alpha_{n-2}\Phi_{n-2}(z))\\
&= (z^2 + z\alpha_{n-2}\bar\beta_1) \Phi_{n-2}(z) -
(z\bar\alpha_{n-2}+\bar\beta_1)\Phi_{n-2}^*(z)\\
&= (z + \alpha_{n-2}\bar\beta_1) (\Phi_{n-1}(z)+
\bar\alpha_{n-2}\Phi_{n-2}^*) - (z\bar\alpha_{n-2}+
\bar\beta_1)\Phi_{n-2}^*(z)\\
&= (z + \alpha_{n-2}\bar\beta_1)\Phi_{n-1}(z) - \bigl\{
\bar\alpha_{n-2}(z+\alpha_{n-2}\bar\beta_1) - (z\bar\alpha_{n-2}+
\bar\beta_1)\bigr\}\Phi_{n-2}^*(z), \end{aligned}
$$
so
\begin{equation} \label{3.7}
\wt\Phi_n(z) = (z+\alpha_{n-2}\bar\beta_1)\Phi_{n-1}(z)-
\bar\beta_1 \rho_{n-2}^2 \Phi_{n-2}^*(z), \quad
\rho_{n-2}^2=1-|\alpha_{n-2}|^2.
\end{equation}
Eliminating $\Phi_{n-2}^*$ from (\ref{3.6}) and (\ref{3.7}) leads
(after elementary computation) to
\begin{equation} \label{3.8}
\Phi_{n-1}(z) = {\wt\Phi_n(z) - A\wt\Phi_{n-1}(z) \over z-B}\,,
\qquad A:= {\bar\beta_1\rho_{n-2}^2 \over \bar\beta_2 -
\bar\alpha_{n-2}} \ne 0,
\end{equation}
with B defined in (\ref{3.5}). It is immediate now from (\ref{3.8})
that if $\zeta_0\in\Sigma(\cC_1)\cap\Sigma(\cC_2)$, then
$\zeta_0=B$.

Let us turn to the intermittency, and begin with (i), which is
referred to as a {\it regular case}. Now $B$ belongs to neither of
$\Sigma(C_l)$, for otherwise it would belong to the other by
(\ref{3.8}), and so $\Sigma(\cC_1)$ and $\Sigma(\cC_2)$ would have
a common point. From (\ref{3.8}) with $z=\zeta_{k,1}$ one has
$$ \Phi_{n-1}(\zeta_{k,1}) = - { A\wt\Phi_{n-1}(\zeta_{k,1})\over
\zeta_{k,1}-B },
$$ and (see (\ref{1.11}))
\begin{equation} \label{3.9}
\mu_{k,1} = {\Phi_{n-1}(\zeta_{k,1}) \over \wt\Phi_n'(\zeta_{k,1})}
= - {A\wt\Phi_{n-1}(\zeta_{k,1}) \over (\zeta_{k,1}-B)
\wt\Phi_n'(\zeta_{k,1}) } >0, \quad k=1,\ldots,n.
\end{equation}
Proposition \ref{propos3} now comes into play with
$$ Z_1=\Sigma(\cC_1)=\{z_{j,1}\}_{j=1}^n; \qquad
Z_2=\Sigma(\cC_2)\cup\{B\}=\{z_{j,2}\}_{j=1}^n,
$$
and $$ z_{j,l}=e^{i x_{j,l}}, \quad j=1,\ldots,n, \ l=1,2; \quad
B=e^{ix_{r,2}},
$$
labeled as in (\ref{3.2}).

Then
$$ \begin{aligned} \mu_{k,1} &= -A\, { {\displaystyle \prod_{j=1}^{n-1}}
(\zeta_{k,1} - \zeta_{j,2}) \over
(\zeta_{k,1}-B) {\displaystyle \prod^{n(k)}}
(\zeta_{k,1}-\zeta_{j,1}) } \\
&= - A\, {{{\displaystyle \prod^{n(r)} e^{i {x_{k,1}+x_{j,2} \over
2} } \prod^{n(r)} 2i \sin {x_{k,1}-x_{j,2} \over 2}} } \over
{\displaystyle e^{i {x_{k,1}+x_{r,2} \over 2} }(2i \sin {x_{k,1} -
x_{r,1} \over 2}) \prod^{n(k)}e^{i {x_{k,1}+x_{j,1} \over 2 } }
\prod^{n(k)} 2i\sin{x_{k,1}-x_{j,1} \over 2} }} \\
&= - {A\over 2i}\, { {\displaystyle \exp \left\{ {i\over2} \sum_1^n
x_{j,2} \right\}} \over {\displaystyle e^{ix_{r,2}}\, \exp \left\{
{i\over2} \sum_1^n x_{j,1} \right\} } }\,\frac1{\displaystyle\sin^2
{x_{k,1}-x_{r,2}\over 2}}\, {{\displaystyle\prod_{j=1}^n \sin
{x_{k,1}-x_{j,2}\over 2}} \over {\displaystyle \prod^{n(k)} {\sin
{x_{k,1}-x_{j,1}\over 2} }}} \\
&= \frac{A_1 \omega_k}{\displaystyle
\sin^2{\frac{x_{k,1}-x_{r,1}}2}}>0, \qquad k=1,2,\ldots,n,
 \end{aligned}
$$
$A_1$ does not depend on $k$, so $\omega_k$ from (\ref{3.4}) have
the same sign, and $Z_1$ and $Z_2$ interlace.

For the {\it singular case} (ii) we proceed in a similar way. Now
$$ Z_1=\Sigma(\cC_1)\backslash\{B\} = \{z_{j,1}\}_{j=1}^{n-1}, \ \
Z_2 =\Sigma(\cC_2) = \{z_{j,2} \}_{j=1}^{n-1}, \ \
B=\zeta_{s,1}=\zeta_{r,2}, $$ so
$$ \{z_{j,1}\}=\{e^{ix_{j,1}}\}=\{\zeta_{j,1}\}_{j\ne s},
\quad \{z_{j,2}\}=\{e^{ix_{j,2}}\}=\{\zeta_{j,2}\}.
$$
For $k=1,\ldots,n, \ k\ne s,$ we have as above
$$ \mu_{k,1} = {A_2\over \displaystyle\sin^2{x_{k,1}-x_{r,2}\over 2} } \cdot
{{\displaystyle\prod_{j=1}^{n-1} \sin {x_{k,1}-x_{j,2}\over 2}}
\over {\displaystyle\prod^{(n-1)(k)} \sin {x_{k,1}-x_{j,1}\over
2}} }
>0,
$$
$A_2$ does not depend on $k$, so by Proposition \ref{propos3}, $Z_1$
and $Z_2$ interlace. \hfill $\square$

\bigskip

{\it Remark.} \ Given CMV matrices (\ref{3.1}), the masses of the
$n$-th spectral measure of $\cC_1$ are given by (\ref{3.9}) in the
regular case. As far as the singular case goes, let $B=\zeta_{n,1}$.
Then (\ref{3.9}) still holds for $k=1,\ldots,n-1$, and (\ref{3.8})
gives
\begin{equation} \label{3.10}
\mu_{n,1} = 1 - A {\wt\Phi_{n-1}'(B) \over \wt\Phi_n'(B)}.
\end{equation}
We need this remark later in Theorem \ref{theor6}.

\bigskip

{\bf Example 2}. Let $\cC=\cC(0,\ldots,0,\beta)$ be a CMV matrix
of order $n$. Then (see Example 1) $\wt\Phi_n(z) = z^n-\bar\beta$,
so
$$\Sigma(\cC) = \left\{\exp\left(i\,\frac{\theta+2\pi j}{n}\right)
\right\}_{j=0}^{n-1}, \quad \beta=e^{-i\theta}, \ \ 0\leq
\theta<2\pi.
$$
Hence for CMV matrices (\ref{3.1}) of the form
$$ \cC_1=\cC(0,\ldots,0;\beta_1), \quad \cC_2=\cC(0,\ldots,0;\beta_2), \qquad
\beta_k=e^{-i\theta_k}, \quad 0\leq\theta_k<2\pi,
$$
of order $n$ and $n-1$, respectively, one has
$$\Sigma(\cC_1) = \left\{\exp\left(\frac{\theta_1+2\pi j}{n}
\right) \right\}_{j=0}^{n-1}, \qquad \Sigma(\cC_2) =
\left\{\exp\left(\frac{\theta_2+2\pi j}{n-1} \right)
\right\}_{j=0}^{n-2}.
$$
Now $B= \beta_2/\beta_1=e^{i(\theta_1-\theta_2)}$. It can be
easily checked that the sets $\Sigma(\cC_1)$ and $\Sigma(\cC_2)$
may have at most one common point $B$, and this singular case
occurs if and only if $\beta_2^n=\beta_1^{n-1}$.

\section{Truncated CMV matrices, inverse problem}
\setcounter{section}{4} \setcounter{equation}{0}

We begin with the regular case. Suppose that two point sets
$Z_1=\{z_{j,1}\}_{j=1}^n$ and $Z_2=\{z_{j,2}\}_{j=1}^{n-1}$ on
$\mathbb{T}$ are given, with no common points, and
\begin{equation} \label{4.1}
z_{j,l}=e^{ix_{j,l}}: \ \ x_{1,1} < x_{1,2} < x_{2,1} < x_{2,2} <
\ldots < x_{n-1,1} < x_{n-1,2} < x_{n,1} < x_{1,1}+2\pi.
\end{equation}

We are aimed at the following result
\begin{theorem}\label{theor5}
Let $\zeta=e^{ix}$ be an arbitrary point on the arc
$(z_{n,1},z_{1,1})$, that is, $x_{n,1}<x<x_{1,1}+2\pi$. Then there
is a unique pair of CMV matrices $(\ref{3.1})$ such that
$$ Z_l=\Sigma(\cC_l), \qquad l=1,2,
$$
and $\zeta=B$ with $B$ defined in $(\ref{3.5})$.
\end{theorem}

{\it Proof.} Write
$$ P_1(z)=\prod_{j=1}^n (z-z_{j,1}), \qquad P_2(z)=
\prod_{j=1}^{n-1} (z-z_{j,2}).
$$
We show first, that for an arbitrary $\zeta \in (z_{n,1,},z_{1,1})$
the nonzero numbers
$$ a_k := {P_2(z_{k,1}) \over (z_{k,1}-\z) P_1'(z_{k,1}) } , \quad
k=1,2,\ldots,n,
$$
have the same argument. Indeed, as above in Section 3
$$ \begin{aligned} a_k &= { {\displaystyle\prod_{j=1}^{n-1} (z_{k,1}-z_{j,2})
\over (z_{k,1}-\z) \displaystyle\prod^{n(k)} (z_{k,1} - z_{j,1})
}} = A_3 b_k, \\
b_k &= {\displaystyle\prod_{j=1}^{n-1} \sin {x_{k,1}-x_{j,2} \over
2} \over  {\displaystyle \sin {x_{k,1}-x \over 2}}
\displaystyle\prod^{n(k)} \sin {x_{k,1}-x_{j,1} \over 2}} ,
\end{aligned} $$ $A_3$ does not depend on $k$. It is clear from
(\ref{4.1}) (cf. Proposition \ref{propos3}) that $b_k<0$ for
$k=1,\ldots,n$, so $a_k$ have the same argument, as claimed.

Hence, there is a unique complex number $v=v(\zeta)$ such that
$$\nu_k:=v a_k>0, \qquad \sum_1^n
\nu_k =1.
$$
A probability measure $d\nu = \sum_{k=1}^n \nu_k \delta(z_{k,1})$,
supported on $Z_1$, appears on the scene, and the unique matrix
$\cC_1(\alpha_0,\ldots,\alpha_{n-2};\beta_1)$ having $d\nu$ as its
$n$-th spectral measure, arises by Theorem~\ref{theor1}. A system
of the \SZ\ polynomials $\Phi_0,\ldots,\Phi_{n-1};\wt\Phi_n$,
associated with $\cC_1$, can be constructed, with $\wt\Phi_n=P_1$,
so $Z_1=\Sigma(\cC_1)$. It remains only to choose $\beta_2$ in an
appropriate way.

Consider the Blaschke product of order 1
\begin{equation} \label{4.2}
b(z) = \beta_1 {z-\alpha_{n-2}\over 1-\bar\alpha_{n-2}z }\,,
\end{equation}
and pick $\beta_2$ as a unique solution of $b(\beta_2)=\zeta$.
Clearly $\zeta=B$ (\ref{3.5}). Now the second CMV matrix $\cC_2$
in (\ref{3.1}) is completely determined, along with its \SZ\
polynomials $\{\Phi_0,\ldots,\Phi_{n-2};\wt\Phi_{n-1}\}$. We wish
to show that $\wt\Phi_{n-1} = P_2$.

By the definition of the spectral measure
(\ref{1.11})--(\ref{1.12}), (\ref{3.9}), and $B=\zeta$, one has
$$ \nu_k =-A {\wt\Phi_{n-1}(z_{k,1}) \over (z_{k,1}-\zeta
)\wt\Phi_n'(z_{k,1})}, \quad k=1,\ldots,n.
$$
But on the other hand
$$ \nu_k = v a_k = v\,{P_2(z_{k,1}) \over (z_{k,1}-\zeta) P_1'(z_{k,1}) }
= v\,{P_2(z_{k,1}) \over (z_{k,1}-\zeta) \wt\Phi_n'(z_{k,1}) },
\quad k=1,\ldots,n,
$$
so $vP_2(z_{k,1})=-A\wt\Phi_{n-1}(z_{k,1})$. The polynomial
$Q=vP_2 + A \wt\Phi_{n-1}$ of degree at most $n-1$ has $n$ roots,
so $Q\equiv 0$, $vP_2=-A\wt\Phi_{n-1}$, and since both $P_2$ and
$\wt\Phi_{n-1}$ are monic, then $P_2=\wt\Phi_{n-1}$, as needed.

To prove uniqueness we provide a procedure of how to restore the
matrices (\ref{3.1}) from two spectra $\sum(\cC_l)$ and the value
$B$. Indeed, $\wt\Phi_n$ and $\wt\Phi_{n-1}$ are determined by
$\sum(\cC_1)$ and $\sum(\cC_2)$, respectively, and the masses of
the $n$-th spectral measure $d\mu_1$ for $\cC_1$ come in a unique
way from (\ref{3.9}) with
$$ A^{-1} = -\sum_{k=1}^n {\wt\Phi_{n-1}(\zeta_{k,1}) \over
(\zeta_{k,1}-B)\wt\Phi_n'(\zeta_{k,1})},
$$
so $\cC_1$ is restored. Finally, $\beta_2$ solves $b(\beta_2)=B$,
$b$ in (\ref{4.2}).

The proof is complete. \hfill $\square$

\bigskip

The singular case can be handled along the same line of reasoning.
Now we have two point sets $Z_1=\{z_{j,1}\}_1^n$ and
$Z_2=\{z_{j,2}\}_1^{n-1}$ on $\mathbb{T}$ with one common point,
$$
z_{j,l} =e^{ix_{j,l}}: \quad
x_{1,1}<x_{1,2}<x_{2,1}<x_{2,2}<\ldots<x_{n-1,1}<x_{n-1,2}=x_{n,1}<x_{1,1}+2\pi.
$$

\begin{theorem}\label{theor6}
In the singular case there are infinitely many pairs of CMV
matrices $(\ref{3.1})$ so that $Z_l=\Sigma(\cC_l)$, $l=1,2$.
\end{theorem}

{\it Proof.} \ As above, we put
$$ P_1=\prod_{j=1}^n (z-z_{j,1}),
\qquad P_2 = \prod_{j=1}^{n-1}(z-z_{j,2}),
$$
and consider the numbers
$$ \begin{aligned} a_k &= {P_2(z_{k,1}) \over (z_{k,1} - z_{n,1}) P_1'(z_{k,1})} =
{ \displaystyle\prod_{j=1}^n (z_{k,1}-z_{j,2}) \over (z_{k,1} -
z_{n,1}) \displaystyle\prod_{j=1}^{n(k)}(z_{k,1}-z_{j,1})}=Rb_k,
\\
R &= {\exp\left\{ {i\over2} \displaystyle\sum_{j=1}^{n-1} x_{j,2}
\right\} \over 2i \exp \left\{ {i\over2} x_{n,1} \right\}\, \exp
\left\{{i\over2}\displaystyle\sum_1^nx_{j,1}\right\}} \ne 0, \\
b_k &= { \displaystyle\prod_{j=1}^{n-1} \sin {x_{k,1} - x_{j,2}
\over 2} \over \displaystyle\sin^2 {x_{k,1} - x_{n,1} \over 2}
\displaystyle\prod^{(n-1)(k)} \sin {x_{k,1} - x_{j,1} \over 2} },
\qquad k=1,2,\ldots,n-1,
\end{aligned} $$
$R$ does not depend on $k$. It follows from the interlacing that
$b_k<0$, \ $k=1,\ldots,n-1$. Put now ($z_{n-1,2}=z_{n,1}$)
$$ a_n := {P_2'(z_{n,1})\over
P_1'(z_{n,1})} = {\displaystyle\prod_{j=1}^{n-2}
(z_{n-1,2}-z_{j,2}) \over \displaystyle\prod_{j=1}^{n-1}
(z_{n,1}-z_{j,1}) } = R\, {\displaystyle\prod_{j=1}^{n-2} \sin
{x_{n-1,2}-x_{j,2}\over 2} \over \displaystyle\prod_{j=1}^{n-1}
\sin {x_{n,1}-x_{j,1}\over 2}}
 = R b_n, \quad b_n>0.
$$
Moreover, let
$$ G(z) = {P_2(z) \over (z-z_{n,1}) P_1(z) } = \sum_{j=1}^n {a_j
\over z-z_{j,1} }\,.
$$
Since $G(z) = O(z^{-2}), \ z\to\infty$, then $\sum_1^n a_j =0$ and
so $\sum_1^n b_j = 0$. Hence there is a complex number $v$ (in
fact, infinitely many such numbers), such that
$$ v a_j >0, \quad j=1,2,\ldots,n-1; \quad v \sum_1^{n-1} a_j<1,
$$
so the numbers $\mu_j = v a_j$, $j=1,\ldots,n-1$, $\mu_n=1+v a_n$
satisfy
$$\mu_j>0, \quad j=1,\ldots,n, \qquad \sum_{j=1}^n \mu_j=1.
$$
Again the measure $d\mu = \sum_1^n \mu_j \delta(z_{j,1})$ comes
in, and by Theorem \ref{theor1} there is a unique CMV matrix
$\cC_1=\cC(\alpha_0,\ldots,\alpha_{n-2};\beta_1)$ with the $n$-th
spectral measure $d\mu$, so $Z_1=\Sigma(\cC_1)$. The sequence of
the \SZ\ polynomials $\{\Phi_0,\ldots,\Phi_{n-1};\wt\Phi_n\}$
arises with $\wt\Phi_n = P_1$.

The choice of $\beta_2$ is the same as in Theorem \ref{theor5}:
$b(\beta_2) = z_{n,1}$, so $B=z_{n,1}$. The second CMV matrix
$\cC_2=\cC(\alpha_0,\ldots,\alpha_{n-3};\beta_2)$ emerges, with
the \SZ\ polynomials $\{\Phi_0,\ldots,\Phi_{n-2};\wt\Phi_{n-1}\}$,
and we want to show that $\wt\Phi_{n-1}=P_2$.

For the masses $\mu_k$ we have by (\ref{3.9}) and (\ref{3.10})
$$ \mu_k = - A {\wt\Phi_{n-1}(z_{k,1}) \over (z_{k,1}-z_{n,1})
\wt\Phi_n'(z_{k,1})}\,, \quad k=1,\ldots,n-1, \quad \mu_n=1 - A
{\wt\Phi_{n-1}'(z_{n,1}) \over \wt\Phi_n'(z_{n,1})\,. }
$$
On the other hand, by the construction
$$ \mu_k = v {P_2(z_{k,1}) \over (z_{k,1}-z_{n,1})
P_1'(z_{k,1})}\,, \quad k=1,\ldots,n-1, \quad \mu_n=1 + v
{P_2'(z_{n,1}) \over P_1'(z_{n,1}) }\,,
$$
so
$$ v P_2(z_{k,1}) = - A\wt\Phi_{n-1}(z_{k,1}), \quad
k=1,\ldots,n-1, \quad v P_2'(z_{n,1}) = - A \wt\Phi_{n-1}'
(z_{n,1}).
$$
Hence the polynomial $\pi=vP_2+A\wt\Phi_{n-1}$ of degree at most
$n-1$ vanishes at $z_{k,1}$, $k=1,\ldots,n-1$, and
$\pi'(z_{n,1})=0$. By the Gauss-Lucas theorem $z_{n,1}$ belongs to
the convex hull of $\{z_{k,1}\}_{k=1}^{n-1}$, that is definitely not
the case, because $|z_{n,1}|=1$ and all $z_{k,1}$ are distinct!
Therefore, $\pi\equiv 0$, $vP_2=-A\wt\Phi_{n-1}$, and in fact
$P_2=\wt\Phi_{n-1}$ as both are monic polynomials. The proof is
complete. \hfill $\square$

\bigskip

The problem arises naturally: whether it is possible (under
additional assumptions) to have a unique solution of the inverse
problem in question. Here is a result of that kind.

\begin{theorem}\label{theor7}
Let $\cC_1^{(l)}=\cC(\a_0^{(l)}, \ldots,
\a_{n-2}^{(l)};\b_1^{(l)})$, $\cC_2^{(l)}=\cC(\a_0^{(l)}, \ldots,
\a_{n-3}^{(l)};\b_2^{(l)})$, \ \  $l=1,2$, be two pairs of
solutions in the singular case. Assume that
$\a_{n-2}^{(1)}=\a_{n-2}^{(2)}$. Then $\cC_j^{(1)}=\cC_j^{(2)}$,
$j=1,2$.
\end{theorem}

{\it Proof}. Since the \SZ\ polynomials $\wt\F_n^{(l)}$,
$\wt\F_{n-1}^{(l)}$ are completely determined by the corresponding
spectra then $\wt\F_n^{(1)}=\wt\F_n^{(2)}$,
$\wt\F_{n-1}^{(1)}=\wt\F_{n-1}^{(2)}$, and
$\b_j^{(1)}=\b_j^{(2)}$, $j=1,2$.

Let us turn to formulae (\ref{3.5}) and (\ref{3.8})
$$
\Phi_{n-1}^{(l)}(z) = {\wt\Phi_n^{(l)}(z) -
A^{(l)}\wt\Phi_{n-1}^{(l)}(z) \over z-B^{(l)}}\,, \qquad A^{(l)}=
{\bar\b_1^{(l)}(\rho_{n-2}^{(l)})^2 \over \bar\b_2^{(l)} -
\bar\a_{n-2}^{(l)}}, \ l=1,2. $$ It follows from the assumption
that $A^{(1)}=A^{(2)}$, $B^{(1)}=B^{(2)}$ and hence
$\F_{n-1}^{(1)}=\F_{n-1}^{(2)}$, which completes the proof. \hfill
$\square$

\end{document}